\documentclass{article}

\usepackage{amssymb,amsthm,amsmath}
\usepackage[all]{xy}

\newtheorem{theorem}{Theorem}[section]
\newtheorem{definition}[theorem]{Definition}
\newtheorem{proposition}[theorem]{Proposition}
\newtheorem{corollary}[theorem]{Corollary}

\newcommand{\mP}{\mathcal{P}}
\newcommand{\mQ}{\mathcal{Q}}
\newcommand{\bX}{\mathbf{X}}

\newenvironment{enumeratereferee}
             {\mbox{ }
              \begin{list}{\arabic{enumi}.}
              {
               \usecounter{enumi}
               \setlength{\itemindent}{\parindent}
               \settowidth{\labelwidth}{\arabic{enumi}.}
               \setlength{\leftmargin}{0pt}
              }
             } {\end{list}}

\title{On completeness in a non-Archimedean setting, via firm reflections}
\author{D. Deses\footnote{Aspirant F.W.O. Vlaanderen} and E. Lowen-Colebunders}

\begin{document}

\maketitle

\begin{abstract}
We develop a completion theory for (general) non-Archimedean spaces
based on the theory on "a categorical concept of completion of objects"
as introduced by G.C.L. Br\"ummer and E. Giuli in  \cite{BG}. Our context is
the construct $\mathbf{NA}_0$ of all Hausdorff non-Archimedean spaces and
uniformly continuous maps and $\mathcal{V}$  is the class of all epimorphic
embeddings in $\mathbf{NA}_0$.

We determine the class ${\bf Inj}\ \mathcal{V}$  of all $\mathcal{V}$-injective
objects and we present an internal characterization as "complete objects".  The
basic tool for this characterization is a notion of small collections that in
some sense preserve the inclusion order on the non-Archimedean structure. We
prove that the full subconstruct $\mathbf{CNA}_0$ consisting of all complete
objects forms a firmly $\mathcal{V}$-reflective subcategory. This means that
every object $X$ in $\mathbf{NA}_0$ has a completion which is a
$\mathcal{V}$-reflection $r_{X}:X\to RX$ into the full
subconstruct $\mathbf{CNA}_0$ of "complete spaces". Moreover this completion is
unique (up to isomorphism) in the sense that, considering $L(\mathbf{CNA}_0)$,
the class of all those morphisms $u: X\to Y$ for which $Ru:RX\to RY$ is an
isomorphism, one has that $\mathcal{V}$  is contained in $L(\mathbf{CNA}_0)$.
In fact one even has $\mathcal{V}=L(\mathbf{CNA}_0)$.

Finally we apply our constructions to the classical case of Hausdorff
non-Archimedean uniform spaces, in that case our completion reduces to the
standard one \cite{R1}, \cite{R2}.

\end{abstract}

{\bf 2000 AMS classification:} 54E15, 54B30, 54D35, 26E30, 18G05.

{\bf Keywords:} completeness, firm reflection, injectives, non-Archimedean
space.

\section{Introduction}

Non-Archimedean uniform spaces were introduced in 1950 by Monna in \cite{Mo}.
They play an important role in non-Archimedean Analysis. First of all
the uniformity of the scalar field $\mathbb{K}$, induced by a non-Archimedean
valuation is itself non-Archimedean, and secondly the non-Archimedean
property is preserved by uniform products and subspaces. In fact every
space obtained from the scalar field by any initial construction, is a
non-Archimedean uniform space. Moreover non-Archimedean uniformities
have proven to provide the right notion of uniformizability for
zero-dimensional topological spaces and they inherit all nice features
from the "exemplary completion theory" for uniform spaces, the meaning of which
is explained below.

Quoting B. Banaschewski from his 1955 paper \cite{Ba}, non-Archimedean
structures "belong to the subfield of General Topology that can be
described by means of equivalence relations".
For Monna \cite{Mo} a non-Archimedean uniformity on $X$ is a uniformity
generated by a filter base of equivalence relations on $X$. In our present
setting, we have to allow for somewhat more general non-Archimedean structures.
Basically we will work with stacks of equivalence relations (see Definition
\ref{def1}). The main reason for passing to the more general setting, is to
include some mathematical structures used in the representation theory of
certain systems. For references see for instance G. Aumann's work on contact
relations or more recent work of B. Ganter and R. Wille on formal concept
analysis \cite{Au}, \cite{GW}. These models combine topological and
lattice-theoretical ideas, and often make use of so called Birkhoff closures
\cite{Bi}, which do not necessarily satisfy the finite additivity of the usual
Kuratowski (i.e. topological) closure. This is in particular the case with
models for physical systems, built on a well defined "lattice of properties of
the system" as developed by Aerts in \cite{A} and by Moore in
\cite{Moo}. There is a natural Birkhoff closure corresponding to this lattice
and recently in \cite{ADVdV} it was shown that in this correspondence, the
lattice of "classical physical properties" gives rise to a zero-dimensional
Birkhoff closure space. The non-Archimedean structures considered in this paper
provide the right notion of uniformizability for these zero-dimensional
closures. Moreover, as in the classical case, they are also stable for initial
constructions.

Our main concern in this note will be with the completion theory for
these (general) non-Archimedean spaces. We will apply a "categorical
concept of completion of objects" as developed by G.C.L. Br\"ummer and
E. Giuli in \cite{BG} to the setting of non-Archimedean spaces. These authors
started from the "exemplary" behavior of the usual completion in the
category $\mathbf{X}=\mathbf{UNIF}_0$ of Hausdorff uniform spaces with
uniformly continuous maps. If $\mathcal{V}$ is the class of all dense
embeddings in $\mathbf{X}$ then, as is well known, every object $X$ in
$\mathbf{X}$ has a completion which is a $\mathcal{V}$-reflection $r_{X}:X \to
RX$ into the full subconstruct $\mathbf{R}$ of "complete spaces". Moreover this
completion is unique (up to isomorphism) in the sense that, considering
$L(\mathbf{R})$, the class of all those morphisms $u:X\to Y$ for which $Ru:
RX \to RY$ is an isomorphism, one has that $\mathcal{V}$ is contained in
$L(\mathbf{R})$. In the case of $\mathbf{X}=\mathbf{UNIF}_0$ one even has
$\mathcal{V}=L(\mathbf{R})$. To describe this exemplary behavior of the
subcategory $\mathbf{R}$ of complete objects in $\mathbf{X}$ the authors of
\cite{BG} used the terminology  "$\mathbf{R}$ is firmly
$\mathcal{V}$-reflective in $\mathbf{X}$ ". The title of our paper refers to
this terminology.

Our context will be the construct $\mathbf{NA}_0$ of all Hausdorff
non-Archimedean spaces and uniformly continuous maps, as introduced in
Definitions \ref{def1} and \ref{def2}. First we determine the class
$\mathcal{V}$ of all epimorphic embeddings in $\mathbf{NA}_0$. Contrary to the
classical case, the epimorphisms are no longer the dense maps described by the
underlying closure of the space. In Theorem \ref{thm1} we characterize the
derived closure operator needed to do the job. From the general results in
\cite{BG} we know that if $\mathbf{NA}_0$ admits a firmly
$\mathcal{V}$-reflective subconstruct $\mathbf{R}$ then $\mathbf{R}$ has to
coincide with the class ${\bf Inj}\ \mathcal{V}$ of all $\mathcal{V}$-injective
objects. In Theorem \ref{thm2} we give an internal characterization of these
injective objects. The basic tool for this characterization is the notion of
small collections (i.e. collections of subsets containing arbitrarily small
sets) that in some sense preserve the inclusion order between the equivalence
relations which determine the structure. These are used to define our "complete
objects" and we prove that the full subconstruct consisting of all complete
objects indeed forms a firmly $\mathcal{V}$-reflective subcategory. We add an
explicit description of the "completion". Finally we apply our constructions to
the classical situation of Hausdorff non-Archimedean uniform spaces, in which
case our completion reduces to the standard one \cite{R1}, \cite{R2}.

For categorical terminology we refer to books such as \cite{AHS}, \cite{HS} or
\cite{P} and for terminology and results on closure operators useful references
are the original papers \cite{DG} and \cite{DGT} or the recent book \cite{DT}.

\section{The construct of non-Archimedean spaces}

In this section we develop the context of our completion theory. We introduce
the construct of Hausdorff non-Archimedean spaces and we pay particular
attention to special morphisms in this construct since in the next paragraph
the class of all epimorphic embeddings will play a key role.

\begin{definition}\label{def1}
A non-Archimedean structure $\mathcal{E}$ on a set $X$ is a stack of
equivalence relations on $X$, i.e. a collection $\mathcal{E}$ of equivalence
relations on $X$ satisfying:
$$E\in \mathcal{E},E\subset E', E' \mbox{ equivalence relation on } X
\Rightarrow E'\in \mathcal{E}$$
The couple $(X,\mathcal{E})$ is called a non-Archimedean space.
\end{definition}

Obviously, instead of working with equivalence relations on $X$ one could
consider partitions of $X$ instead. So an alternative for Definition \ref{def1}
is to consider a set $\beta$ of partitions of $X$ satisfying:
{\it $$\mP\in \beta,\mP\prec \mP', \mP' \mbox{ partition on } X \Rightarrow
\mP'\in \beta$$}
where $\prec$ is the refinement relation defined on covers $\mP\prec
\mathcal{P'}$ ($\mP$ refines $\mP'$) iff $\forall P\in\mP:\exists
P'\in\mP':P\subset P'$.

We will use the following notations. A  non-Archimedean space
will be written as $\bX$ and we shall use $\beta_{X}$ and $\mathcal{E}_{X}$ to
refer to the corresponding structures and $X$ for the underlying set. We will
write $[x]_\mP$ as well as $E[x]$ for the equivalence class of a point $x\in
X$.

Remark that every partition star refines itself. Therefore every
non-Archi\-me\-dean space provides a base for a pre-nearness space in the
sense of \cite{P}. It is even a base for a uniform semi-nearness space as in
\cite{DLCVdV}. Moreover if $\mathcal{E}_{X}$ is closed under finite intersections
then it forms a base for a collection of entourages of a non-Archimedean
uniform space \cite{Ba}, \cite{Mo}, \cite{R1}, \cite{R2}. In this case
$\beta_{X}$ is closed under the operation $\wedge$ given by $\mP\wedge
\mQ=\{P\cap Q|P\in\mP,Q\in \mQ\}$ and it generates a uniform nearness space in
the sense of \cite{H}.

\begin{definition}\label{def2}
A function $f:\bX\to \mathbf{Y}$ between non-Archimedean spaces is
called uniformly continuous if
$$\forall E\in \mathcal{E}_{Y}:(f\times f)^{-1}(E)\in \mathcal{E}_{X}$$
\end{definition}

In terms of partitions this is obviously equivalent to
$$\forall \mP\in\beta_{Y}:f^{-1}(\mP)\in \beta_{X}$$

The category of non-Archimedean spaces together with the uniformly continuous
maps will be denoted $\mathbf{NA}$. It is a topological construct in the sense
of \cite{AHS}. This means that initial (and final) structures for arbitrary
class indexed sources (and sinks) can be formed in $\mathbf{NA}$. In particular
the objects on a fixed underlying set $X$, form a complete lattice with largest
element the discrete object $\mathbf{D}_{X}$ and smallest one the indiscrete
object $\mathbf{I}_{X}$.

Given a source $(f_i:X\to \bX_i)_{i\in I}$ in $\mathbf{NA}$ then the initial
structure $\mathcal{E}$ is given by  $\{(f_i\times f_i)^{-1}(E_i)|i\in I,
E_i\in\mathcal{E}_{X_i}\}$.

$\mathbf{NA}_0$ is the subconstruct consisting of the $T_0$ objects of
$\mathbf{NA}$. Applying the usual definition \cite{Ma} we say that $\bX$ is a
$T_0$ object if and only if every uniformly continuous map from
the indiscrete object $\mathbf{I}_{\{0,1\}}$ to $\mathbf{X}$ is constant. This
equivalently means that for any two different points $x$ and $y$ in $X$ there
is an equivalence relation $E\in\mathcal{E}_{X}$ such that $E[x]\not =E[y]$.

In view of this separation condition the objects in $\mathbf{NA}_0$ will be
called Hausdorff non-Archimedean spaces. From the results of Marny in
\cite{Ma} it follows that $\mathbf{NA}_0$ is an extremally epireflective
subconstruct of $\mathbf{NA}$ and as such it is initially structured in the
sense of \cite{N}, \cite{P}. In particular $\mathbf{NA}_0$ is complete,
cocomplete and well powered, it is an (epi - extremal mono) category and an
(extremal epi - mono) category \cite{HS}. Also from the general setting
\cite{N} it follows that the monomorphisms in $\mathbf{NA}_0$ are exactly the
injective uniformly continuous maps and a morphism in $\mathbf{NA}_0$ is an
extremal epimorphism if and only if it is a regular epimorphism if and only if
it is surjective and final.

In order to describe the epimorphisms and the extremal monomorphisms in
$\mathbf{NA}_0$ the following result on cogenerators in $\mathbf{NA}_0$ is very
useful.

\begin{proposition}\label{prop1}
Let $\bX$ be a Hausdorff non-Archimedean space. We have that
$$i:\bX\to \Pi_{\mP\in \beta_{X}} \mathbf{D}_\mP: x\to([x]_\mP)_{\mP\in \beta_{X}}$$
is an embedding (i.e. injective and initial).
\end{proposition}

\begin{proof}
Consider the source:
$$(i_\mP:X\to\mathbf{D}_\mP:x \to[x]_\mP)_{\mP\in \beta_{X}}$$
Since $(i_\mP\times i_\mP)^{-1}(\Delta_\mP)=E_\mP$, where $E_\mP$ denotes
the equivalence relation defined by $\mP$, we know that $\beta_{X}$ is the
initial structure for this source. Because $\bX$ is Hausdorff we have that the
above source is point separating. Hence
$$i:\bX\to \Pi_{\mP\in \beta_{X}} \mathbf{D}_\mP: x\to([x]_\mP)_{\mP\in \beta_{X}}$$
is an embedding.
\end{proof}

This means that Hausdorff non-Archimedean spaces form exactly the epireflective
hull in $\mathbf{NA}$ of the class of all discrete objects, or using the
terminology of \cite{BGH}, $\mathbf{NA}_0$ is cogenerated, with respect to all
embeddings, by the class of all discrete objects, i.e. the Hausdorff
non-Archimedean spaces are exactly the subspaces of a product of discrete
spaces.

To simplify notations we shall write $\Pi \beta_{X}$ instead of $\Pi_{\mP\in
\beta_{X}} \mathbf{D}_\mP$.

As we will see next, as opposed to the classical case, the epimorphisms in
$\mathbf{NA}_0$ can not be described as the "dense" uniformly continuous maps,
with "denseness" defined by the underlying closure. We need to determine the
$\mathbf{NA}_0$-regular closure operator as introduced in \cite{DG}, \cite{DGT}
and define "denseness" accordingly.

\begin{definition}\label{def3} \cite{DG}, \cite{DGT}
Given a non-Archimedean space $\mathbf{X}$ and a subset \mbox{$M\subset X$}, a
point $x$ of $X$ is in the regular closure of $M$ in $\mathbf{X}$ iff \\
for every Hausdorff non-Archimedean space $\mathbf{Z}$ and for every pair of
uniformly continuous maps $f,g:\mathbf{X}\to\mathbf{Z}$
$$f_{|M}=g_{|M}\Rightarrow f(x)=g(x)$$
in this case we write $x\in reg_{X}(M)$.
\end{definition}

Using Proposition \ref{prop1} we immediately obtain the following equivalent
description of Definition \ref{def3}.\\
{\it $x\in reg_{X}(M)$ iff for every discrete space $\mathbf{D}$ and for every
pair of uniformly continuous maps $f,g:\mathbf{X}\to\mathbf{D}$
$$f_{|M}=g_{|M}\Rightarrow f(x)=g(x)$$}

In order to obtain an explicit description of the regular closure operator we
introduce the following notation.

\begin{definition}\label{def4}
Let $\bX$ be a non-Archimedean space and $M\subset X$.\\
$x\in \zeta_{X}(M)$ iff for every two equivalence relations $E_1,E_2\in
\mathcal{E}_{X}$, which coincide on $M$, we have that $E_1[x]\cap E_2[x]\cap M
\not =\emptyset$.
\end{definition}

\begin{theorem}\label{thm1}
For every non-Archimedean space $\mathbf{X}$ and $M\subset X$:
$$\zeta_{X}(M)=reg_{X}(M)$$
\end{theorem}

\begin{proof}
Let $\bX$ be a non-Archimedean space, let $x\in X$, $M\subset X$ such that
$x\not\in reg_{X}(M)$. There is a discrete object $\mathbf{D}$ and there
are uniformly continuous maps $f,g:\bX\to \mathbf{D}$ for which
$f_{|M}=g_{|M}$ and $f(x)\not = g(x)$. Consider \mbox{$E_1=(f\times
f)^{-1}(\Delta_D),E_2=(g\times g)^{-1}(\Delta_D)$}. Clearly $E_1$ and
$E_2$ belong to $\mathcal{E}_{X}$ but do not satisfy the condition in Definition
\ref{def4}

\begin{sloppy}
Conversely, if $x\not \in \zeta_{X}(M)$, choose $E_1,E_2\in \mathcal{E}_{X}$
such that $E_1$ and $E_2$ coincide on $M$ and  $E_1[x]\cap E_2[x]\cap
M=\emptyset$. Let \mbox{$C=\{E_1[m]\cap M|m\in M\}$} \mbox{$=\{E_2[m]\cap
M|m\in M\}$} and $D=C \cup \{a,b\}$ where $a,b\not \in C$. Write $\mathbf{D}$
for the discrete object on $D$. We define the following functions:
$$f:\bX\to \mathbf{D}:y \mapsto \Bigl \{
\begin{array}{cc}
E_1[m]\cap M & \mbox{if } \exists m \in M: (y,m)\in E_1 \\
a & \mbox {if } \forall m \in M: (y,m)\not \in E_1 \\
\end{array}$$
and
$$g:\bX\to \mathbf{D}:y \mapsto \Bigl \{
\begin{array}{cc}
E_2[m]\cap M & \mbox{if } \exists m \in M: (y,m)\in E_2 \\
b & \mbox {if } \forall m \in M: (y,m)\not \in E_2 \\
\end{array}$$
Clearly $f$ and $g$ are uniformly continuous, they coincide on $M$ but
\mbox{$f(x)\not =g(x)$}.
\end{sloppy}
\end{proof}

It follows from \cite{DGT} that
$$\zeta:\{\zeta_{X}:\mathcal{P}(X)\to\mathcal{P}(X)\}_{\mathbf{X}\in|\mathbf{NA}|}$$
defines a closure operator on $\mathbf{NA}$.

Clearly $\zeta$ is hereditary in the sense that for a space $\mathbf{Y}$ and a
subspace $\mathbf{X}$ in $\mathbf{NA}$ and $M\subset X\subset Y$, we have
$\zeta_{X}(M)=\zeta_{Y}(M)\cap X$.

A subset $M$ of a non-Archimedean space $\mathbf{X}$ is $\zeta$-dense
($\zeta$-closed) if \mbox{$\zeta_{X}(M)=X$} ($\zeta_{X}(M)=M$). A map
$f:\mathbf{X}\to\mathbf{Y}$ between non-Archimedean spaces is $\zeta$-dense
($\zeta$-closed) if $f(X)$ has the corresponding property with respact to $Y$
\cite{DG}, \cite{DGT}.

\begin{proposition}\label{prop4}
In $\mathbf{NA}_0$ we have:
\begin{enumerate}
\item The epimorphisms are exactly the $\zeta$-dense uniformly continuous
maps.
\item The extremal monomorphisms coincide with the regular monomorphisms and
they both coincide with the $\zeta$-closed embeddings.
\end{enumerate}
\end{proposition}

\begin{proof}
\begin{enumeratereferee}

\item This follows from Theorem 2.8 in \cite{DGT} since $\zeta$ is the regular
closure operator determined by $\mathbf{NA}_0$.

\item Using proposition 2.6 in \cite{DGT} and the fact that $\mathbf{NA}_0$ is
extremally epireflective in $\mathbf{NA}$, for an $\mathbf{NA}_0$ morphism
$f:\mathbf{X}\to\mathbf{Y}$ the following implications hold.
(i) $f$ $\zeta$-closed embedding $\Rightarrow$ (ii) $f$ regular monomorphism
$\Rightarrow$ (iii) $f$ extremal monomorphism $\Rightarrow$ (iv) $f$ embedding.
To see that (iii) implies that $f$ is $\zeta$ closed, we have to use (weak)
hereditariness of $\zeta$, either by applying 6.2 in \cite{DT} or by the
following direct argument. Let $M=\zeta_{Y}(f(X))$ and
$h:\mathbf{M}\to\mathbf{Y}$ the associated $\zeta$-closed embedding. Then there
exists a unique map $g$ such that
$$\xymatrix{
\mathbf{X}\ar[rd]_g \ar[rr]^f & &\mathbf{Y} \\
& \mathbf{M} \ar[ur]_h & \\
}$$
By the (weak) hereditariness of the $\zeta$-closure, $g$ is $\zeta$-dense and
so it is an epimorphism in $\mathbf{NA}_0$. It follows that $g$ is an
isomorphism and then $f$ is $\zeta$-closed.
\end{enumeratereferee}

\end{proof}

From the previous characterization of the epimorphisms in $\mathbf{NA}_0$ it
now follows that $\mathbf{NA}_0$ is cowell-powered.

It suffices to observe that given an $\zeta$-dense map $f:\mathbf{X}\to
\mathbf{Y}$, with $\mathbf{X}$ fixed, there is a one to one correspondence
between $\beta_{f(X)}$ and $\beta_{Y}$. Therefore the cardinality of $Y$ is
uniformly bounded.

\section{The firm $\mathcal{V}$-reflective subconstruct consisting of complete
objects of $\mathbf{NA}_0$}

In this paragraph we show that $\mathbf{NA}_0$ admits a completion theory that
is "exemplary" in the sense explained in the introduction.

Let $\mathcal{V}$ be the class of epimorphic embeddings of $\mathbf{NA}_0$, by
Proposition \ref{prop4} $\mathcal{V}$ consists of all $\zeta$-dense embeddings.
This class $\mathcal{V}$ satisfies the following conditions
\begin{itemize}
\item[($\alpha$)] closedness under composition
\item[($\beta$)] closedness under composition with isomorphisms on both sides
\end{itemize}
($\alpha$) and ($\beta$) are standing assumptions made in \cite{BG} and enable
us to apply to $\mathbf{NA}_0$ the theory developed in that paper.
We will prove that $\mathbf{NA}_0$ admits a $\mathcal{V}$-reflective
subconstruct $\mathbf{R}$ which is firm in the terminology of \cite{BG}.
Explicitly this means that:
\begin{enumerate}
\item Every object $\bX$ has a reflection $r_{X}:\bX \to R\bX$ into $\mathbf{R}$
such that $r_{X}\in \mathcal{V}$.
\item If $L(\mathbf{R})$ is the class of all morphisms $u:\bX \to \mathbf{Y}$
in $\mathbf{NA}_0$ for which $Ru:R\mathbf{X}\to R\mathbf{Y}$ is an isomorphism
then $\mathcal{V}=L(\mathbf{R})$.
\end{enumerate}

By proposition 1.6 in \cite{BGH}, in order to construct $\mathbf{R}$ it
suffices to find a class consisting of $\mathcal{V}$-injective objects that
cogenerates $\mathbf{NA}_0$, with respect to embeddings (cf. Proposition
\ref{prop1}).

Recall that a Hausdorff non-Archimedean space $\mathbf{B}$ is
$\mathcal{V}$-injective if for each $v:\bX\to \mathbf{Y}$ in $\mathcal{V}$ and
$f:\bX \to \mathbf{B}$ uniformly continuous there exists a uniformly continuous
$f':\mathbf{Y} \to \mathbf{B}$ such that $f'\circ v=f$. In this case $f'$ is
called an extension of $f$ along $v$. ${\bf Inj} \ \mathcal{V}$ denotes the full
subcategory of all $\mathcal{V}$-injective objects in $\mathbf{NA}_0$.

\begin{proposition}\label{prop5}
Every discrete Hausdorff non-Archimedean space is
$\mathcal{V}$-injective.
\end{proposition}

\begin{proof}
Let $u:\bX\to\mathbf{Y}$ be a $\zeta$-dense embedding between two Hausdorff
non-Archimedean spaces, and let $f:\bX\to\mathbf{D}$ be a uniformly continuous
function to a discrete space. By the initiality of $u$ we have an $E\in
\mathcal{E}_{Y}$ such that $(u\times u)^{-1}(E)=(f\times
f)^{-1}(\Delta_D)$. For $y\in Y$ we choose $x_y\in X$ such that $(y,u(x_y))\in
E$. We define $$f^*:\mathbf{Y}\to\mathbf{D}:y\mapsto f(x_y)$$
Clearly $f^*$ is a well-defined uniformly continuous map which is an extension
of $f$ along $u$.
\end{proof}

Combining Propositions \ref{prop1} and \ref{prop5} and using Theorem 1.6 in
\cite{BGH} and Theorem 1.4 and Proposition 1.14 in \cite{BG} we can formulate
the next result:

\begin{proposition} \label{prop6}
The following hold:
\begin{enumerate}
\item $\mathbf{NA}_0$ admits a unique firmly $\mathcal{V}$-reflective
subconstruct $\mathbf{R}$.
\item $\mathbf{R}$ coincides with the class ${\bf Inj}\ \mathcal{V}$.
\item $\mathbf{R}$ coincides with the epireflective hull in $\mathbf{NA}_0$ of
the class of all discrete spaces.
\end{enumerate}
\end{proposition}

We next present an internal characterization of the objects in the firm
$\mathcal{V}$-reflective subconstruct. In order to do this we formulate the
following.

\begin{definition}
Let $\bX$ be a non-Archimedean space. As usual, a choice function is a map
$f:\beta_{X}\to \cup \beta_{X}$ such that for any $\mP\in \beta_{X}$ one has that
$f(\mP)\in \mP$. A choice function is order preserving iff $\mP\prec \mQ$
implies $f(\mP)\subset f(\mQ)$.
\end{definition}

\begin{definition}
A Hausdorff non-Archimedean space is complete iff
for every order preserving choice function $f$ there is an $x \in
\cap_{\mP\in \beta_{X}}f(\mP)$.
\end{definition}

The point $x$ is called a limit point of $f$ and we will say that $f$ converges
to $x$. Note that in this case the limit point $x$ is unique.

The following proposition links this concept of completeness to the firmly
$\mathcal{V}$-reflective subcategory we described before.

\begin{proposition}\label{prop2}
Let $\bX$ be a Hausdorff non-Archimedean space and consider the
embedding $i:\bX \to \Pi \beta_{X}$ from Proposition \ref{prop1}. We have the
following equivalence.
$$z=(z_\mP)_{\mP\in \beta_{X}}\in \zeta_{\Pi \beta_{X}}(i(X)) \iff
\forall\mP,\mQ\in\beta_{X}: \mP\prec \mQ \Rightarrow z_\mP\subset z_\mQ$$
\end{proposition}

\begin{proof}
First let $z\in \zeta_{\Pi\beta_{X}}(i(X))$ and let $\mP\prec\mQ$, where
$\mP,\mQ\in\beta_{X}$. For \mbox{$E_\mP=(pr_\mP\times
pr_\mP)^{-1}(\Delta_\mP)$ and $E_\mQ=(pr_\mQ\times
pr_\mQ)^{-1}(\Delta_\mQ)$}, we know that there is an $x\in X$ such that
$i(x)\in E_\mP[z]\cap E_\mQ[z]$. Therefore $[x]_\mP=z_\mP\in\mP$ and
$[x]_\mQ=z_\mQ\in\mQ$. Because $\mP\prec\mQ$ we know that there is a $Q\in \mQ$
such that $z_\mP\subset Q$. Since $z_\mQ\in \mQ$ and both $Q$ and $z_\mQ$
contain $x$ we have that $z_\mQ=Q$. Hence $z_\mP\subset z_\mQ$.

Conversely, choose $E_1,E_2$ in $\mathcal{E}_{\Pi\beta_{X}}$ such that both
coincide on $i(X)$. Clearly one can write $E_1=(pr_\mP\times
pr_\mP)^{-1}(E_\mP)$ and $E_2=(pr_\mQ\times pr_\mQ)^{-1}(E_\mQ)$ where
$E_\mP$ and $E_\mQ$ are equivalences on $\mP$ and $\mQ$ respectively.

For the equivalence relation $E_\mP$, consider classes
$[x]_\mP,[y]_\mP\in \mP$ of points $x$ and $y$ in $X$, we have
$$\begin{array}{ccc}
([x]_\mP,[y]_\mP)\in E_\mP &\iff &(i(x),i(y))\in E_1 \\
&\iff &(i(x),i(y))\in E_2 \\
&\iff &([x]_\mQ,[y]_\mQ)\in E_\mQ\\
\end{array}$$
Hence $\mathcal{R}=\{\{y\in X| ([x]_\mP,[y]_\mP)\in E_\mP\}|x\in X\}$ is a
partition of $X$, such that $\mP,\mQ\prec \mathcal{R}$. So by the
hypothesis we have $z_\mP,z_\mQ\subset z_\mathcal{R}$.
Thus there is an $x\in X$ for which $(i(x),z)\in E_1$ and
$(i(x),z)\in E_2$. Finally we have $z\in \zeta_{\Pi\beta_{X}}(i(X))$.

\end{proof}

\begin{theorem}
Let $\bX$ be a non-Archimedean space. The following are equivalent:
\begin{enumerate}
\item $\bX$ is complete.
\item $\bX$ is a $\zeta$-closed subspace of $\Pi\beta_{X}$.
\item $\bX$ is $\mathcal{V}$-injective.
\end{enumerate}
\end{theorem}

\begin{proof} We prove the following implications:\\

\noindent $1\Rightarrow 2:$ Suppose $\bX$ is complete. We prove that $i(X)$ is
$\zeta$-closed in $\Pi\beta_{X}$. Let \mbox{$z=(z_\mP)_{\mP\in \beta_{X}}\in
\zeta_{\Pi\beta_{X}}(i(X))$}, by Proposition \ref{prop2} we know that
\mbox{$f:\beta_{X}\to\cup\beta_{X}:$} $\mP\mapsto z_\mP$ is an order preserving
choice function. Clearly the completeness of $\bX$ guarantees the existence of
$x\in X$ such that $z=i(x)$.\\

\noindent $2\Rightarrow 3:$ Let $\bX$ be a $\zeta$-closed subspace of $\Pi
\beta_{X}$. Then since $\mathbf{NA}_0$ is complete, well-powered and
cowell-powered, applying 37.6 in \cite{HS}: $\bX$ belongs to the epireflective
hull of all discrete spaces. By Proposition \ref{prop6} $\bX$ is
$\mathcal{V}$-injective.\\

\noindent $3\Rightarrow 1:$ Suppose $\bX$ is $\mathcal{V}$-injective in
$\mathbf{NA}_0$. In view of Proposition \ref{prop4} (2) we can conclude that
$X$ is $\zeta$-closed in every Hausdorff non-Archimedean space in which it is
embedded. In particular $i:\bX \to \Pi \beta_{X}$ is a $\zeta$-closed embedding.
Now if $f$ is any order preserving choice function, Proposition \ref{prop2}
implies that $f$ converges to some point of $X$.

\end{proof}

Let $\mathbf{CNA}_0$ be the full subconstruct of $\mathbf{NA}_0$ consisting of
the complete Hausdorff non-Archimedean spaces. $\mathbf{CNA}_0$ is the unique
$\mathcal{V}$-reflective subconstruct of $\mathbf{NA}$. In the next proposition
we give an explicit description of the reflection, for a space $\bX$ in
$\mathbf{NA}_0$, which is in fact the "unique completion" of $\bX$.

\begin{theorem}\label{thm2}
Let $\bX$ be a Hausdorff non-Archimedean space. Let $\hat{X}$ be
the set consisting of all order preserving choice functions of $\bX$. \\
For every $E\in \mathcal{E}_{X}$, let $\mP$ be the partition given by $E$. We
define $\hat{E}=\{(f,g)\in \hat{X}|f(\mP)=g(\mP)\}$. Consider the
non-Archimedean space $\hat{\bX}$ where the structure is given by the
equivalence relations on $\hat{X}$ that contain a relation of
\mbox{$\hat{\mathcal{E}}=\{\hat{E}| E\in\mathcal{E}_{X}\}$}.  We have the
following: \begin{enumerate} \item $\hat{\bX}$ is a complete Hausdorff
non-Archimedean space. \item $\bX$ is a $\zeta$-dense subspace of $\hat{\bX}$.
\end{enumerate} \end{theorem}

\begin{proof}
\begin{enumeratereferee}
\item Suppose that $f,g\in \hat{X}$ are different. Then there is a $\mP\in
\beta_{X}$ such that $f(\mP)\not =g(\mP)$. So $\hat{\bX}$ is Hausdorff.

Let $\hat{f}$ be an order preserving choice function on $\hat{\bX}$. We make an
order preserving choice function on $\bX$ as follows. Each $\mP\in \beta_{X}$ has
a corresponding equivalence relation $E_\mP$, for which $\hat{E_\mP}$ has a
partition $\hat{\mP}$ on $\hat{X}$. We define $f:\beta_{X}\to \cup
\beta_{X}:\mP\mapsto \hat{f}(\hat{\mP})$. For any $\hat{\mP}$ we have that
$f\in \hat{f}(\hat{\mP})$. Hence $\hat{f}$ converges to $f$, so $\hat{\bX}$ is
complete.

\item Consider the following map:
$$j:\bX\to \hat{\bX}:x\to (f_x:\beta_{X}\to \cup \beta_{X}:\mP \mapsto
[x]_\mP)$$
Clearly this $f_x$ always is an order preserving choice function. Since $\bX$
is Hausdorff $j$ obviously is injective.

Let $E\in \mathcal{E}_{X}$, clearly $E=(j\times j)^{-1}(\hat{E})$. Therefore
$j$ is initial.

Let $f\in \hat{X}$ and let $F_1,F_2\in \mathcal{E}_{\hat{X}}$ which
coincide on $j(X)$. There exist $E_1,E_2\in \mathcal{E}_{X}$ and the
corresponding partitions $\mP_1,\mP_2$ such that $\hat{E}_1\subset F_1$ and
$\hat{E}_2\subset F_2$.

We have that $\hat{E}_1[f]=\{g\in \hat{X}|f(\mP_1)=g(\mP_1)\}$. Since
\mbox{$\emptyset \not = f(\mP_1)\in \mP_1$} there is an $x\in f(\mP_1)$ such
that $f_x\in \hat{E}_1[f]$, so $(f,f_x)\in \hat{E}_1\subset F_1$. Analogously
there is a $y\in f(\mP_2)$ for which $(f,f_y)\in \hat{E}_2 \subset F_2$.

Since $F_1,F_2$ coincide on $j(X)$, we have that $(j\times j)^{-1}(F_1)=
(j\times j)^{-1}(F_2)$, the latter corresponding to a $\mP'\in \beta_{X}$
for which $\mP_1,\mP_2\prec \mP'$. Since $f$ is order preserving and because of
our choice of $x$ and $y$ we have that $x,y\in f(\mP')$. Hence $(f_x,f_y)\in
F_2$ and then also $(f,f_x)\in F_2$. Therefore $f_x\in F_1[f]\cap F_2[f]\cap
j(X)$. So finally $f\in \zeta_{\hat{X}}(j(X))$.

\end{enumeratereferee}
\end{proof}

From the previous theorem and starting with $\bX$ in $\mathbf{NA}_0$ we now can
conclude that with respect to the class $\mathcal{V}$ of $\zeta$-dense
embeddings, $(\hat{X},j)$ is the unique completion of $\bX$. Indeed if
$r_{X}:\bX \to R\bX$ is the reflection of $\bX$ in $\mathbf{CNA}_0$, then with
the notations of Theorem \ref{thm2} we can consider the diagram:
$$\xymatrix{
R\bX \ar[dr]^{r(j)}& & \\
\bX \ar[u]_{r_{X}} \ar[r]_j &\hat{\bX}\simeq R\hat{\bX} \\
}$$
Since $\mathbf{CNA}_0$ is firmly $\mathcal{V}$-reflective and $j$ is
$\zeta$-dense, this means that $r(j)$ is an isomorphism.

From this it follows that a uniformly continuous map $u:\bX \to \mathbf{Y}$
from a Hausdorff non-Archimedean space $\bX$ to a complete Hausdorff
non-Archimedean space $\mathbf{Y}$ can be uniquely extended to a uniformly
continuous map $\hat{u}:\hat{\bX}\to \mathbf{Y}$. We describe this extension
explicitly as follows.

Let $f\in \hat{X}$. For any $\mP\in \beta_{Y}$ we know that $u^{-1}(\mP)\in
\beta_{X}$, hence $f(u^{-1}(\mP))$ is a class of $u^{-1}(\mP)$ and
$u(f(u^{-1}(\mP)))$ is a subset of a class from $\mP$, which we will write as
$f_u(\mP)$. This defines an order preserving choice function as follows:
$$f_u:\beta_{Y} \to \cup \beta_{Y}:\mP\mapsto f_u(\mP)$$
The extension of $u$ is then given by:
$$\hat{u}:\hat{\bX}\to \mathbf{Y}:f\mapsto lim \ f_u$$
where $lim \ f_u$ is the unique limit of $f_u$, which exists since
$\mathbf{Y}$ is Hausdorff and complete.

\section{The case of non-Archimedean uniform spaces}

In this last section we will now show that the previously described completion
is in fact a generalization of the classical completion of a Hausdorff
non-Archimedean uniform space, as described in \cite{Ba},\cite{R1},\cite{R2}.

Let $\mathbf{X}$ be a Hausdorff non-Archimedean uniform space as
introduced by Monna in \cite{Mo}, described by a collection of entourages
$\mathcal{U}$. We write $\mathcal{E}$ for the collection of
equivalence relations in $\mathcal{U}$ and $\beta$ for the corresponding
collection of partitions. Obviously $\mathbf{X}$ uniquely corresponds
to a Hausdorff non-Archimedean space in our sense.

\begin{proposition}\label{prop7}
Let $\bX$ be a non-Archimedean uniform space and let $\mathcal{F}$ be a minimal
Cauchy filter in $\mathbf{X}$. Then there is a unique element $(z_\mP)_{\mP\in
\beta}\in \Pi \beta$ such that $\mathcal{F}=stack \ \{z_\mP|\mP\in \beta\}$.
\end{proposition}

\begin{proof}
We only have to check uniqueness since such a collection exists because $\mathcal{F}$ is
minimal. If $\mathcal{F}=stack \ \{z_\mP|\mP\in \beta\}=stack \ \{z'_\mP|\mP\in
\beta\}$, then for every $\mP$ $z_\mP\cap z'_\mP$ is nonempty since
$\mathcal{F}$ is a filter. Hence $z_\mP=z'_\mP$.
\end{proof}

\begin{proposition}
Let $\bX$ be a non-Archimedean uniform space. Then there is a one to one
correspondence between the order preserving choice functions of $\bX$ and the
minimal Cauchy filters of $\bX$.
\end{proposition}

\begin{proof}
Let $A$ be the set of all order preserving choice functions of $\bX$ and let
$B$ denote the set of all its minimal Cauchy filters. The following maps
describe the needed one to one correspondence.
$$F:A \to B:f\mapsto \mathcal{F}_f$$
where $\mathcal{F}_f=stack \ \{z_\mP|\mP\in\beta_{X}\}$ with $z_\mP=f(\mP)$ for
$\mP\in\beta_{X}$.
$$G:B\to A:\mathcal{F}\mapsto f_\mathcal{F}$$
where $f_\mathcal{F}(\mP)=z_\mP$ is uniquely defined by Proposition
\ref{prop7}.

By definition $\mathcal{F}_f$ contains arbitrarily small sets. For $\mP,\mQ\in
\beta_{X}$ we have $z_{\mP\wedge\mQ}\subset z_\mP\cap z_\mQ$, so
$\mathcal{F}_f$ is a filter. By the same argument as in the proof of
Proposition \ref{prop7} one has that $\mathcal{F}_f$ is minimal.

$f_\mathcal{F}$ is a well defined choice function by Proposition \ref{prop7}.
For $\mP,\mQ\in \beta_{X}$ we have $z_{\mP\wedge\mQ}=z_\mP\cap z_\mQ$. It follows
that $f_\mathcal{F}$ is order preserving.

After a simple verification one sees that $F$ and $G$ are bijective and
inverse to one another.

\end{proof}

Since through the bijections $F$ and $G$ convergent order preserving choice
functions correspond to convergent minimal Cauchy filters we can conclude the
following.

\begin{corollary}
The completion as developed in Theorem \ref{thm2}, when applied to a Hausdorff
non-Archimedean uniform space reduces to the standard completion.
\end{corollary}


\begin{thebibliography}{}

\bibitem{AHS}
{\scriptsize J. Ad\'amek, H. Herrlich and G. Strecker,
\emph{Abstract and Concrete Categories},
Wiley, New York (1990).}

\bibitem{A}
{\scriptsize D. Aerts,
\emph{Foundations of quantum physics: a general realistic and operational approach},
Internat. J. Theoret. Phys. \textbf{38(1)} (1999), 289--358.}

\bibitem{ADVdV}
{\scriptsize D. Aerts, D. Deses and A. Van der Voorde,
\emph{Connectedness applied to closure spaces and state property systems},
J. Elec. Eng. \textbf{52} (2001), 18--21.}

\bibitem{Au}
{\scriptsize G. Aumann,
\emph{Kontaktrelationen},
Bayer. Akad. Wiss. Math.-Nat. Kl. Sitzungsber. (1970), 67--77.}

\bibitem{Ba}
{\scriptsize B. Banaschewski,
\emph{\"Uber nulldimensionale R\"aume},
Math. Nachr., \textbf{13} (1955), 129--140.}

\bibitem{Bi}
{\scriptsize G. Birkhoff,
\emph{Lattice Theory},
American Mathematical Society, Providence, Rhode Island (1940).}

\bibitem{BG}
{\scriptsize G. C. L. Br{\"u}mmer and E. Giuli,
\emph{A categorical concept of completion of objects},
Comment. Math. Univ. Carolin. \textbf{33(1)} (1992), 131--147.}

\bibitem{BGH}
{\scriptsize G. C. L. Br{\"u}mmer, E. Giuli and H. Herrlich,
\emph{Epireflections which are completions},
Cahiers Topologie G\'eom. Diff. Cat\'eg. \textbf{33(1)} (1992), 71--73.}

\bibitem{DLCVdV}
{\scriptsize D. Deses, E. Lowen-Colebunders, A. De Groot-Van der Voorde,
\emph{Extensions of closure spaces},
submitted for publication.}

\bibitem{DG}
{\scriptsize D. Dikranjan and E. Giuli,
\emph{Closure operators. I.}
Topology Appl. \textbf{27} (1987), 129--143.}

\bibitem{DGT}
\scriptsize{D. Dikranjan, E. Giuli and A. Tozzi,
\emph{Topological categories and closure operators},
Quaestiones Math. \textbf{11(3)} (1988), 323--337.}

\bibitem{DT}
\scriptsize{D. Dikranjan and W. Tholen,
\emph{Categorical structure of closure operators},
Kluwer Academic Publ. Dordrecht \textbf{346} (1995).}

\bibitem{GW}
{\scriptsize B. Ganter and R. Wille,
\emph{Formal Concept Analysis},
Springer Verlag, Berlin (1998).}

\bibitem{H}
{\scriptsize H. Herrlich,
\emph{Topological structures},
Math. Centre Tracts, \textbf{52} (1974), 59-122.}

\bibitem{HS}
{\scriptsize H. Herrlich and G. Strecker,
\emph{Category theory},
Sigma Series in Pure Mathematics, Heldermann Verlag, Berlin (1987).}

\bibitem{Ma}
{\scriptsize Th. Marny,
\emph{On epireflective subcategories of topological categories},
General Topology Appl. \textbf{10(2)} (1979), 175--181.}

\bibitem{Mo}
{\scriptsize A. F. Monna,
\emph{Remarques sur les m\'etriques non-Archim\'ediennes I \& II},
Indag. Math. \textbf{12} (1950), 122--133 \& 179--191.}

\bibitem{Moo}
{\scriptsize D. J. Moore,
\emph{Categories of representations of physical systems},
Helv. Phys. Acta \textbf{68} (1995), 658--678.}

\bibitem{N}
{\scriptsize L. D. Nel,
\emph{Initially structured categories and Cartesian closedness},
Canad. J. Math. \textbf{27(6)} (1975), 1361--1377.}

\bibitem{P}
{\scriptsize G. Preuss,
\emph{Theory of Topological Structures},
D. Reidel Publishing Company, Dordrecht (1988).}

\bibitem{R1}
{\scriptsize A. C. M. van Rooij,
\emph{Non-Archimedean uniformities},
Kyungpook Math. J., \textbf{10} (1970), 21-30.}

\bibitem{R2}
{\scriptsize A. C. M. van Rooij,
\emph{Non-Archimedean functional analysis},
Pure and Appl. Math., Marcel Decker Inc. New York (1987).}

\end{thebibliography}
\end{document}